\documentclass[11pt,reqno]{amsart}
\usepackage{amssymb,verbatim,enumerate,ifthen}
\usepackage[mathscr]{eucal}
\usepackage[utf8]{inputenc}
\usepackage[T1]{fontenc}
\def\N{\mathbb{N}}
\def\R{\mathbb{R}}
\def\Q{\mathbb{Q}}
\def\Z{\mathbb{Z}}
\def\F{\mathbb{F}}

\newtheorem{theorem}{Theorem}
\newtheorem*{theorem*}{Theorem}
\def\Thm#1#2{\ifthenelse{\equal{#1}{*}}{\begin{theorem*}#2\end{theorem*}}
             {\begin{theorem}\label{T#1}#2\end{theorem}}}
\newtheorem{Atheorem}{Theorem}

\def\THM#1#2{\begin{Atheorem}\label{T#1}#2\end{Atheorem}}
\def\thm#1{Theorem~\ref{T#1}}

\newtheorem{proposition}[theorem]{Proposition}
\newtheorem*{proposition*}{Proposition}
\def\Prp#1#2{\ifthenelse{\equal{#1}{*}}{\begin{proposition*}#2\end{proposition*}}
             {\begin{proposition}\label{P#1}#2\end{proposition}}}

\newtheorem{corollary}[theorem]{Corollary}
\newtheorem*{corollary*}{Corollary}
\def\Cor#1#2{\ifthenelse{\equal{#1}{*}}{\begin{corollary*}#2\end{corollary*}}
             {\begin{corollary}\label{C#1}#2\end{corollary}}}
\newtheorem{Acorollary}[Atheorem]{Corollary}

\def\COR#1#2{\begin{Acorollary}\label{C#1}#2\end{Acorollary}}
\def\cor#1{Corollary~\ref{C#1}}

\newtheorem{lemma}[theorem]{Lemma}
\newtheorem*{lemma*}{Lemma}
\def\Lem#1#2{\ifthenelse{\equal{#1}{*}}{\begin{lemma*}#2\end{lemma*}}
             {\begin{lemma}\label{L#1}#2\end{lemma}}}
\newtheorem{Alemma}[Atheorem]{Lemma}

\def\lem#1{Lemma~\ref{L#1}}

\theoremstyle{definition}
\newtheorem{remark}[theorem]{Remark}
\newtheorem*{remark*}{Remark}
\def\Rem#1#2{\ifthenelse{\equal{#1}{*}}{\begin{remark*}\rm #2\end{remark*}}
             {\begin{remark}\label{R#1}\rm #2\end{remark}}}

\newtheorem{example}[theorem]{Example}
\newtheorem*{example*}{Example}
\def\Exa#1#2{\ifthenelse{\equal{#1}{*}}{\begin{example*}\rm #2\end{example*}}
             {\begin{example}\label{Ex#1}\rm #2\end{example}}}

\def\eq#1{{\rm(\ref{E#1})}}
\def\Eq#1#2{\ifthenelse{\equal{#1}{*}}
  {\begin{equation*}\begin{aligned}#2\end{aligned}\end{equation*}}
  {\begin{equation}\begin{aligned}\label{E#1}#2\end{aligned}\end{equation}}}

\begin{document}
\begin{flushright}
\textit{Submitted to: Aequationes Math.}
\end{flushright}
\vspace{5mm}

\date{\today}

\title[Characterizations of higher-order convexity]
{Characterizations of higher-order convexity properties with respect to Chebyshev systems}

\author[Zs. P\'ales]{Zsolt P\'ales}
\author[\'E. Sz\'ekelyn\'e Rad\'acsi]{\'Eva Sz\'ekelyn\'e Rad\'acsi}
\address{Institute of Mathematics, University of Debrecen,
H-4010 Debrecen, Pf.\ 12, Hungary}
\email{\{pales,radacsi.eva\}@science.unideb.hu}


\thanks{This research has been supported by the Hungarian Scientific Research Fund (OTKA) Grant K111651}

\dedicatory{Dedicated to the 70th birthday of Professor Roman Ger}

\begin{abstract}
In this paper various notions of convexity of real functions with respect to Chebyshev
systems defined over arbitrary subsets of the real line are introduced. As an auxiliary notion, a concept of
a relevant divided difference and also a related lower Dinghas type derivative are also defined. The main
results of the paper offer various characterizations of the convexity notions in terms of the nonnegativity of
a generalized divided difference and the corresponding lower Dinghas type derivative.
\end{abstract}

\maketitle

\section{Introduction}

Denote by $\N$, $\Z$, $\Q$, $\R$ the sets of natural, integer, rational, and real numbers, respectively.
Given a set $H\subseteq\R$, the set of positive elements of $H$ is denoted by $H_+$. Thus, for instance,
$\N=\Z_+$.

For a set $H\subseteq\R$, denote the simplex of strictly ordered $n$-tuples of the elements of $H$ by
$\sigma_n(H)$, i.e.,
\Eq{*}{
  \sigma_n(H):=\{(x_1,\dots,x_n)\in H^n\mid x_1<\dots<x_n\}.
}
Obviously, $\sigma_n(H)\neq\emptyset$ if and only if the cardinality $|H|$ of $H$ is at least $n$.

Provided that $|H|\geq n$, for a vector valued function $ \omega=(\omega_1, \dots, \omega_n): H \to \R^n $,
the functional operator $\Phi_\omega:\sigma_n(H)\to\R$ is defined by
\Eq{*}{
    \Phi_\omega(x_1,\dots,x_n):=
   \left|\begin{array}{ccc}
     \omega_1(x_1) &  \dots & \omega_1(x_n)\\
     \vdots & \ddots &  \vdots \\
     \omega_n(x_1) &  \dots & \omega_n(x_n)
   \end{array}\right|
   \qquad\big((x_1,\dots,x_n)\in\sigma_n(H)\big).
}
We say that $ \omega $ is an \emph{$n$-dimensional positive (resp.\ negative) Chebyshev system over $H$} if
$\Phi_\omega$ is positive (resp.\ negative) over $\sigma_n(H)$, respectively.

The following systems are the most important particular cases for positive Chebyshev systems. For more
important examples we refer to the books by Karlin \cite{Kar68} and Karlin--Studden \cite{KarStu66}.
\begin{enumerate}[(i)]
\item The function $\omega:\R\to\R^n$ given by $\omega(x):=(1, x, \dots, x^{n-1}) $ is an
$n$-dimensional positive Chebyshev system on $ \R $. Indeed, then $\Phi_\omega(x_1,\dots,x_n)$ is the well-known
Vandermonde determinant, i.e.,
    \Eq{*}{
        \Phi_\omega(x_1,\dots,x_n)=
        \left|\begin{array}{ccc}
            1 &  \dots & 1\\
            x_1  & \dots & x_n\\
            \vdots & \ddots & \vdots \\
            x_1^{n-1} & \dots & x_n^{n-1}
        \end{array}\right|
        = \prod_{1\leq i < j \leq n}{(x_j-x_i)}
        \qquad \big((x_1,\dots,x_n) \in \sigma_n(\R)\big),
    }
    and this determinant is obviously positive on $\sigma_n(\R)$. The above Chebyshev system is called the
\emph{standard, or polynomial $n$-dimensional Chebyshev system}.
\item The function $ \omega(x):=\big( 1,\cos(x), \sin(x), \dots, \cos(nx), \sin(nx)\big) $ is a
$(2n+1)$-dimensional positive Chebyshev system on any open interval $I$ whose length is less than or equal to $ 2\pi $.
Indeed, after some calculation, for $ (x_1,\dots,x_{2n+1}) \in \sigma_{2n+1}(I) $, we get
    \Eq{*}{
        \Phi_\omega (x_1,\dots,x_{2n+1})
        =4^{n^2}\prod_{1\leq k<j\leq 2n+1}{\sin\Big(\frac{x_j-x_k}{2}\Big)}.
    }
    Provided that the length of the open interval $I$ is
less than or equal to $2\pi$, for $1\leq k<j\leq 2n+1$, we have $0<\frac{x_j-x_k}{2}<\pi$, hence
$\sin\big(\frac{x_j-x_k}{2}\big)>0$. Therefore, $\Phi_\omega (x_1,\dots,x_{2n+1})>0$.
\item The function $ \omega(x):=\big( \cos(x), \sin(x), \dots, \cos(nx), \sin(nx)\big) $
is a $(2n)$-dimensional positive Chebyshev system on any open interval $I$ whose length is less than or equal
to $ \pi $. Indeed, after some calculation, for $ (x_1,\dots,x_{2n}) \in \sigma_{2n}(I) $, we get
    \Eq{*}{
        \Phi_\omega (x_1,\dots,x_{2n})
        =\frac{4^{n(n-1)}}{(n!)^2}\sum_{(j_1,\dots,j_{2n})}
        \prod_{k=1}^n \cos\Big(\frac{x_{j_{k+1}}-x_{j_k}}{2}\Big)
        \prod_{1\leq k<j\leq 2n}{\sin\Big(\frac{x_j-x_k}{2}\Big)},
    }
    where the summation is taken over all permutations
$(j_1,\dots,j_{2n}) $ of the set $\{1,\dots,2n\}$. Provided that the length of the open interval $I$ is less
than or equal to $\pi$, for $1\leq k<j\leq 2n$, we have $0<\frac{x_j-x_k}{2}<\frac\pi2$, hence
$\sin\big(\frac{x_j-x_k}{2}\big)>0$ and $\cos\big(\frac{x_j-x_k}{2}\big)>0$. Therefore, $\Phi_\omega
(x_1,\dots,x_{2n})>0$.
\item For the function $ \omega(x):=(1,x^2 ) $, we get that $\Phi_\omega(x_1,x_2)=(x_1+x_2)(x_2-x_1)$.
Therefore $\omega$ is a $2$-dimensional positive Chebyshev system on $ \R_+ $, but it is not a Chebyshev
system on $ \R $ (observe that $ \Phi_{\omega}(-1,1)=0$).
\end{enumerate}

Given a positive Chebyshev system $\omega:H\to\R^n$, a function $f:H\to\R$ is called
\emph{$\omega$-convex} (i.e., \emph{convex with respect to the Chebyshev system $\omega$}) if,
$\Phi_{(\omega,f)}$ is nonnegative over $\sigma_{n+1}(H)$. A function $f:H\to\R$ is strictly $\omega$-convex
if a function $ (\omega_1,\dots,\omega_n,f)$ is an $(n+1)$-dimensional positive Chebyshev system over $H$.

For $k\geq0$, define the $k$th power function $p_k:\R\to\R$ by $p_k(x):=x^k$. As we have seen it before,
$(p_0,\dots,p_{n-1})$ is an $n$-dimensional Chebyshev system. The notion of convexity with respect to this
system, called polynomial convexity, was introduced by Hopf \cite{Hop26} and by Popoviciu
\cite{Pop44}. The particular case, when $\omega=(p_0,p_1)$, simplifies to the notion of standard convexity,
moreover, for $(x,y,z)\in\sigma_3(H)$, the inequality $\Phi_{(p_0,p_1,f)}(x,y,z)\geq0$ is equivalent to
\Eq{*}{
 f(y)\leq \frac{z-y}{z-x}f(x)+\frac{y-x}{z-x}f(z).
}

The next result is well-known for the standard convexity (\cite{NikPal03}), it provides the main motivation
for our paper.

\THM{1a}{
Let $x_0<\dots<x_m$ be real numbers, $f:\{x_0,\dots,x_m\}\to\R$, where $m\geq2$. If, for all
$i\in\{1,\dots,m-1\}$, the function $f$ is convex on the 3-element set $\{x_{i-1},x_i,x_{i+1}\}$, i.e., $f$
fulfills the inequality
    \Eq{*}{
        \Phi_{(p_0,p_1,f)}(x_{i-1},x_i,x_{i+1})\geq 0
        \qquad \big(i\in\{1,\dots,m-1\}\big),
    }
then $f$ is convex over $\{x_0,\dots,x_m\}$, i.e., for every $0\leq k<\ell<n\leq m$,
    \Eq{*}{
        \Phi_{(p_0,p_1,f)}(x_k,x_\ell,x_n)\geq 0.
    }
}

In order to derive the proof of this result the notion of second-order divided differences, its connection to
convexity and the Chain Inequality established in \cite{NikPal03}. Recall, that the second-order divided
difference of a function $f:H\to\R$ (where $|H|\geq3$), is defined by
\Eq{*}{
  \big[x,y,z;f\big]
   :=\frac{f(x)}{(y-x)(z-x)}+\frac{f(y)}{(x-y)(z-y)}+\frac{f(z)}{(x-z)(y-z)}
      \qquad\big((x,y,z)\in\sigma_3(H)\big).
}
It is easy to see that
\Eq{2dd}{
  \big[x,y,z;f\big]=\frac{\Phi_{(p_0,p_1,f)}(x,y,z)}{\Phi_{(p_0,p_1,p_2)}(x,y,z)}
  \qquad\big((x,y,z)\in\sigma_3(H)\big).
}
The Vandermonde determinant $\Phi_{(p_0,p_1,p_2)}(x,y,z)$ being positive, the inequality
$\Phi_{(p_0,p_1,f)}(x,y,z)\geq0$ is equivalent to $\big[x,y,z;f\big]\geq0$. On the other hand, in terms of
second-order divided differences, the decomposition holds true.

\THM{A}{
Let $x_0<\dots<x_m$ be real numbers, where $m\geq2$. Then, for every $0\leq k<\ell<n\leq m$, there exist
$A_i\in [0,1]$ with $\sum_{i=1}^{m-1}A_i=1$ such that, for all function $f:\{x_0,\dots,x_m\}\to\R$,
\Eq{*}{
  \big[x_{k},x_\ell,x_{n};f\big]
  =\sum_{i=1}^{m-1}A_i\big[x_{i-1},x_i,x_{i+1};f\big].
}}

This theorem easily implies

\COR{A}{{\rm (Chain Inequality)}
Let $x_0<\dots<x_m$ be real numbers, $f:\{x_0,\dots,x_m\}\to\R$, where $m\geq2$. Then, for every
$0\leq k<\ell<n\leq m$,
\Eq{*}{
  \min_{i\in\{1,\dots,m-1\}}\big[x_{i-1},x_i,x_{i+1};f\big]
  \leq \big[x_{k},x_\ell,x_{n};f\big]
  \leq \max_{i\in\{1,\dots,m-1\}}\big[x_{i-1},x_i,x_{i+1};f\big].
}}

Applying the above corollary and formula \eq{2dd}, one can easily obtain a proof for \thm{1a}.

The above \thm{A} can be extended to higher-order divided differences applying Lemma 2.6 in Chapter XV of
the book \cite{Kuc85} by M. Kuczma. Using this extension, generalizations of \cor{A} and \thm{1a} can easily
be derived to the case when the underlying Chebyshev system is $(p_0,\dots,p_{n-1})$. Based on these results,
A. Gil\'anyi and Zs. P\'ales \cite{GilPal08} introduced various generalizations of Popoviciu's convexity
notions and obtained several characterizations of them.

In what follows, in Section 2 we shall introduce a notion of divided difference which is suitable for the
applications in the setting of a general Chebyshev system $\omega$. Then we prove a result which will be
analogous to \thm{A} and has consequences that generalize \cor{A} and \thm{A}. In Section 3 we introduce the
notions of $(t,\omega)$-convexity and $\omega$-Jensen convexity and establish several implications among
these convexity properties. Finally, in Section 4, we introduce various lower Dinghas type derivatives and
obtain mean value inequalities in three different settings for them. These results immediately yield the
characterizations of convexity properties in terms of the nonnegativity of the relevant lower Dinghas type
derivative and show also that these convexity properties are localizable ones. The results obtained in this
paper directly generalize those in \cite{NikPal03} and in \cite{GilPal08}.

\section{Convexity over discrete sets}

Throughout this section, let $n\in\N$, and let $H\subseteq\R$ with $|H|\geq n+2$ and let
$\omega=(\omega_1,\dots, \omega_n):H\to\R$ be an $n$-dimensional positive Chebyshev system and choose
$\omega_{n+1}:H\to\R$ such that $\omega_{n+1}$ is strictly convex with respect to $\omega$, that is,
$\omega^*=(\omega_1,\dots,\omega_n,\omega_{n+1})$ is an $(n+1)$-dimensional positive Chebyshev system.

For a function $f:H\to\R$, the generalized \emph{$\omega^*$-divided difference $[x_0,\dots,x_n;f]_{\omega^*}$
of $f$} is defined by
    \Eq{*}{
        \big[x_0,\dots,x_n;f\big]_{\omega^*}
        :=\frac{\Phi_{(\omega,f)}(x_0,\dots,x_n)}{\Phi_{\omega^*}(x_0,\dots,x_n)}
             \qquad \big( (x_0,\dots,x_n) \in \sigma_{n+1}(H) \big).
    }
Clearly, if $\omega^*=(\omega_1,\dots, \omega_n,\omega_{n+1})=(p_0,\dots,p_{n-1},p_n)$, then,
in view of the identity \eq{2dd}, $\big[x_0,\dots,x_n;f\big]_{\omega^*}$ is equal to the standard
$n$th-order divided difference $\big[x_0,\dots,x_n;f\big]$. On the other hand, given a positive Chebyshev
system $\omega=(\omega_1,\dots, \omega_n)$, there is no unique candidate for the strictly
$\omega$-convex function $\omega_{n+1}$.

For the proof the main theorem below, we need the following auxiliary result.

\Lem{1}{Let $ x_0<\dots<x_{n+1} $ be arbitrary elements of $H$. Then, for all $0\leq j \leq n+1$, there exists
a constant $A_j \in [0,1]$ such that, for every function $f:\{x_0,\dots,x_{n+1} \}\to\R$, the following
equality holds
\Eq{L1}{
        \big[x_0,\dots,x_{j-1},x_{j+1},\dots,x_{n+1};f\big]_{\omega^*}
         =A_j\big[x_0,\dots,x_n;f\big]_{\omega^*} + (1-A_j)\big[x_1,\dots,x_{n+1};f\big]_{\omega^*}.
}}

\begin{proof}
If $j=0$ or $j=n+1$ the statement of the lemma trivially holds with $A_0:=0$ and $A_{n+1}:=1$, respectively.
Let $0< j < n+1$ be arbitrary and define $A_j$ by the equality
    \Eq{Aj}{
        \frac{A_j}{1-A_j}:=
        \frac{\Phi_{\omega^*}(x_0,\dots,x_n)}{\Phi_{\omega^*}(x_1,\dots,x_{n+1})}\cdot
        \frac{\Phi_{\omega}(x_1,\dots,x_{j-1},x_{j+1},\dots,x_{n+1})}
             {\Phi_{\omega}(x_0,\dots,x_{j-1},x_{j+1},\dots,
x_n)}.
    }
Since the right side of this equality is positive, hence $\frac{A_j}{1-A_j}>0$, which yields $A_j \in ]0,1[$.

To prove that \eq{L1} holds for all $f:\{x_0,\dots,x_{n+1} \}\to\R$, using that $\omega^*$ is an
$(n+1)$-dimensional Chebyshev system, we can find constants $\alpha_1,\dots,\alpha_{n+1} \in \R$ such
that
    \Eq{*}{
        (\alpha_1\omega_1+\cdots+\alpha_{n+1}\omega_{n+1})(x_i)=f(x_i) \qquad  i\in \{ 0,\dots,n+1 \} \setminus
        \{j\}.
    }
Then, with the notation $g:=f- (\alpha_1\omega_1+\cdots+\alpha_{n+1}\omega_{n+1})$, we can decompose
$f$ into the following form
    \Eq{*}{
        f=(\alpha_1\omega_1+\cdots+\alpha_{n+1}\omega_{n+1}) + g,
    }
where $ g(x_i)=f(x_i)-(\alpha_1\omega_1+\cdots+\alpha_{n+1}\omega_{n+1})(x_i) = 0$
for $ i\in \{ 0,\dots,n+1 \} \setminus \{j\}$. Substituting $f$ into the left hand side of the equation
\eq{L1}, by obvious transformations rules of determinants, we get
    \Eq{*}{
        &\big[x_0,\dots,x_{j-1},x_{j+1},\dots,x_{n+1};f\big]_{\omega^*}\\
        &=\big[x_0,\dots,x_{j-1},x_{j+1},\dots,x_{n+1};
              \alpha_1\omega_1+\cdots+\alpha_{n+1}\omega_{n+1}\big]_{\omega^*}
        +\big[x_0,\dots,x_{j-1},x_{j+1},\dots,x_{n+1};g\big]_{\omega^*}
        =\alpha_{n+1}.
    }
On the other hand, using \eq{Aj}, the right hand side of the equation \eq{L1} reduces to the following form
    \Eq{*}{
        &A_j\big[x_0,\dots,x_n;f\big]_{\omega^*} + (1-A_j)\big[x_1,\dots,x_{n+1};f\big]_{\omega^*}\\
        &=A_j \bigg( \alpha_{n+1}+ (-1)^{n+j+2}g(x_j)
        \frac{\Phi_{\omega}(x_0,\dots,x_{j-1},x_{j+1},\dots,x_n)}
             {\Phi_{\omega^*}(x_0,\dots,x_n)} \bigg)\\
        &\qquad+ (1-A_j) \bigg(  \alpha_{n+1}+(-1)^{n+j+1}g(x_j)
        \frac{\Phi_{\omega}(x_1,\dots,x_{j-1},x_{j+1},\dots,x_{n+1})}
             {\Phi_{\omega^*}(x_1,\dots,x_{n+1})}\bigg) \\
        &=\alpha_{n+1}+(-1)^{n+j}g(x_j)
          \bigg(A_j\frac{\Phi_{\omega}(x_0,\dots,x_{j-1},x_{j+1},\dots,x_n)}
                  {\Phi_{\omega^*}(x_0,\dots,x_n)}
                -(1-A_j)\frac{\Phi_{\omega}(x_1,\dots,x_{j-1},x_{j+1},\dots,x_{n+1})}
             {\Phi_{\omega^*}(x_1,\dots,x_{n+1})}\bigg) \\
        &=\alpha_{n+1}.
    }
This completes the proof of \eq{L1}.
\end{proof}

The first main theorem of our paper extends an analogous result known for standard higher-order divided
differences, cf.\ \cite{Kuc85}.

\Thm{1}{Let $m\in\N$ with $n\leq m\leq |H|-1$ and let $ x_0<\dots<x_{m} $ be elements of $H$. Then, for all $0
\leq i_0<\dots<i_n\leq m$, there exist $A_i \in [0,1]$ with $\sum_{i=0}^{m-n}A_i=1$ such that, for all
function $f:\{x_0,\dots,x_m\}\to\R$,
\Eq{1a}{
    \big[x_{i_0},\dots,x_{i_n};f\big]_{\omega^*}
     =\sum_{i=0}^{m-n}A_i\big[x_i,\dots,x_{i+n};f\big]_{\omega^*}
}
holds.}

\begin{proof}
The above theorem can be proved by induction with respect to $m-n$. For $m-n=0$,
the theorem trivially holds with $A_0=1$. If $m-n=1$, then the theorem results from \lem{1}.

Now assume the theorem is true for a $k=m-n$. To prove that the theorem is valid for $k+1=m-n$, let $
x_0<\dots<x_{n+k+1} $ be arbitrary elements of $H$.

Take an element $ x_{i_{n+1}}$ of $H$, which is different from $ x_{i_0},\dots,x_{i_n} $. Let
$ \{x_{i_0}^{'},\dots,x_{i_{n+1}}^{'} \} $ be a set such that,
$\{x_{i_0}^{'},\dots,x_{i_{n+1}}^{'}\}=\{x_{i_0},\dots,x_{i_{n+1}}\}$ and
$x_{i_0}^{'}<\dots<x_{i_{n+1}}^{'}$.
By \lem{1}, there exist constants $B_0,B_1 \in [0,1] $ with $B_0+B_1=1$ such that, for all function
$f:\{x_{i_0}^{'},\dots,x_{i_{n+1}}^{'}\}\to\R$,
\Eq{*}{
    \big[x_{i_0},\dots,x_{i_n};f\big]_{\omega^*}=
    B_0 \big[x_{i_0}^{'},\dots,x_{i_n}^{'};f\big]_{\omega^*}
    + B_1 \big[x_{i_1}^{'},\dots,x_{i_{n+1}}^{'};f\big]_{\omega^*}.
}
Notice that $ \{x_{i_0}^{'},\dots,x_{i_n}^{'}\}$ is a subset of $
\{x_0,\dots,x_{n+k}\} $ and $ \{x_{i_1}^{'},\dots,x_{i_{n+1}}^{'}\}$ is a subset of $ \{x_1,\dots,x_{n+k+1}\}$
and in both case we have $n+k+1-(n+1)=k$. By the induction hypothesis there exist constants $ C_i,D_i \in
[0,1] $ such that, $ \sum_{i=0}^{k}C_i=\sum_{i=1}^{k+1}D_i=1$ and for all function
$f:\{x_{0},\dots,x_{n+k+1}\}\to\R$,
\Eq{*}{
    \big[x_{i_0}^{'},\dots,x_{i_n}^{'};f\big]_{\omega^*}
    =\sum_{i=0}^{k}C_i\big[x_i,\dots,x_{i+n};f\big]_{\omega^*},
    \\
    \big[x_{i_1}^{'},\dots,x_{i_{n+1}}^{'};f\big]_{\omega^*}
    =\sum_{i=1}^{k+1}D_i\big[x_i,\dots,x_{i+n};f\big]_{\omega^*}.
}
Hence,
\Eq{*}{
    &\big[x_{i_0},\dots,x_{i_n};f\big]_{\omega^*}=
    B_0 \bigg(\sum_{i=0}^{k}C_i\big[x_i,\dots,x_{i+n};f\big]_{\omega^*}\bigg)
    + B_1 \bigg(\sum_{i=1}^{k+1}D_i\big[x_i,\dots,x_{i+n};f\big]_{\omega^*}\bigg)\\&=
    B_0C_0\big[x_0,\dots,x_{n};f\big]_{\omega^*}+\sum_{i=1}^{k}(B_0C_i+B_1D_i)\big[x_i,\dots,x_{i+n};f\big]_{\omega^*} + B_1D_{k+1}\big[x_{k+1},\dots,x_{n+k+1};f\big]_{\omega^*}.
}
Finally, we have to show that the sum of the coefficient is equal to $1$.
\Eq{*}{
    B_0C_0+\sum_{i=1}^{k}(B_0C_i+B_1D_i)+ B_1D_{k+1}= B_0\sum_{i=0}^{k}C_i+B_1\sum_{i=1}^{k+1}D_i=B_0+B_1=1.
}
Thus we have obtained \eq{1a} for $k+1$. Induction ends the proof.
\end{proof}

As an immediate consequence, we obtain

\Cor{1}{(Generalized Chain Inequality)
Let $m\in\N$ with $n\leq m\leq |H|-1$ and let $ x_0<\dots<x_{m} $ be elements of $H$. Then, for
all $0 \leq i_0<\dots<i_n\leq m$ and for all function $f:\{x_0,\dots,x_m\}\to\R$,
\Eq{C1}{
    \min_{0\leq i\leq m-n}\big[x_i,\dots,x_{i+n};f\big]_{\omega^*}
    \leq\big[x_{i_0},\dots,x_{i_n};f\big]_{\omega^*}
    \leq\max_{0\leq i\leq m-n}\big[x_i,\dots,x_{i+n};f\big]_{\omega^*}.
}}

\begin{proof} By \thm{1}, there exists $A_i \in [0,1]$ with $\sum_{i=0}^{m-n}A_i=1$ such that \eq{1a} holds
for all functions $f:\{x_0,\dots,x_m\}\to\R$. On the other hand, it is easy to see that
\Eq{*}{
  \min_{0\leq i\leq m-n}\big[x_i,\dots,x_{i+n};f\big]_{\omega^*}
    \leq\sum_{i=0}^{m-n}A_i\big[x_i,\dots,x_{i+n};f\big]_{\omega^*}
    \leq\max_{0\leq i\leq m-n}\big[x_i,\dots,x_{i+n};f\big]_{\omega^*}.
}
These inequalities combined with \eq{1a} yield \eq{C1}.
\end{proof}

\Cor{2}{Let $m\in\N$ with $n\leq m\leq |H|-1$, let $ x_0<\dots<x_{m} $ be elements of $H$, and
let $f:\{x_0,\dots,x_m\}\to\R$. If, for all $0\leq i \leq m-n$, $f$ is $\omega$-convex on the $(n+1)$-element
set $\{x_i,\dots,x_{i+n}\}$, i.e.,
\Eq{*}{
    \Phi_{(\omega,f)}(x_i,\dots,x_{i+n})\geq 0,
}
then $f$ is $\omega$-convex on $\{x_0,\dots,x_{m}\}$, i.e., for all
$ 0 \leq i_0 <\dots<i_n \leq m$, the following inequality holds
  \Eq{*}{
    \Phi_{(\omega,f)}(x_{i_0},\dots,x_{i_n}) \geq 0.
  }
}

\begin{proof} Choose $\omega_{n+1}:H\to\R$ such that
$\omega^*=(\omega_1,\dots,\omega_n,\omega_{n+1})$ be a positive Chebyshev system over $H$. Since
$ \Phi_{(\omega,f)}(x_i,\dots,x_{i+n})\geq 0 $, hence $[x_i,\dots,x_{i+n};f]_{\omega^*}\geq 0$ for all
$ i\in\{0,\dots,m-n\}$. According to the \thm{1} (or to \cor{1}), this implies that
$[x_{i_0},\dots,x_{i_n};f]_{\omega^*}\geq 0 $ for all $ 0 \leq i_0 <\dots<i_n \leq m$.
Hence, $ \Phi_{(\omega,f)}(x_{i_0},\dots,x_{i_n}) \geq 0 $ holds for all $ 0 \leq i_0 <\dots<i_n \leq m$.
\end{proof}

\section{$(t,\omega)$-convexity and $\omega$-Jensen convexity}

Let $I\subset\R$ be a nondegenerate interval, $n\in\N$, and let $\omega:I\to\R^n$ be an $n$-dimensional
positive Chebyshev system over $I$ throughout this section.

For $t= (t_1,\dots,t_n) \in \R_+^n$ and for a permutation $\pi$ of the index set $\{1,\dots,n\}$, define
$t\circ\pi$ by $t\circ\pi:=(t_{\pi(1)},\dots,t_{\pi(n)})$. A function $f:I\to\R$ is said to be
\emph{$(t,\omega)$-convex on $I$} if
\Eq{to}{
    \Phi_{(\omega,f)}(x,x+t_1h,\dots,x+(t_1+\dots+t_n)h)\geq 0
}
holds for all $h>0$, $x\in I$ with $x+(t_1+\dots+t_n)h\in I$. We call $f$ \emph{cyclically
$(t,\omega)$-convex on $I$} if it is $(t\circ\pi,\omega)$-convex for all cyclic permutations $\pi$ of
$\{1,\dots,n\}$. We call $f$ \emph{symmetrically $(t,\omega)$-convex on $I$} if it is
$(t\circ\pi,\omega)$-convex for all permutations $\pi$ of $\{1,\dots,n\}$. Finally, we call $f$
\emph{$\omega$-Jensen convex on $I$}, if it is $(t,\omega)$-convex with $t=\mathbf{1}_n$, where 
$\mathbf{1}_k$ stands for the vector $(1,\dots,1)\in\R^k$ for $k\in\N$.
Note that, for $n\geq2$, with the substitution $y:=x+(t_1+\dots+t_n)h$,
\eq{to} is satisfied for all $h>0$, $x\in I$ with $x+(t_1+\dots+t_n)h\in I$ if and only if
\Eq{*}{
  \Phi_{(\omega,f)}\Big(x,\tfrac{t_2+\dots+t_n}{t_1+\dots+t_n}x+\tfrac{t_1}{t_1+\dots+t_n}y,
       \dots,\tfrac{t_n}{t_1+\dots+t_n}x+\tfrac{t_1+\dots+t_{n-1}}{t_1+\dots+t_n}y,y\Big)\geq 0
}
holds for all $x,y\in I$ with $x<y$. In particular, $f$ is $\omega$-Jensen convex on $I$ if and only if
\Eq{*}{
  \Phi_{(\omega,f)}\big(x,\tfrac{n-1}{n}x+\tfrac{1}{n}y,\dots,\tfrac{1}{n}x+\tfrac{n-1}{n}y,y\big)\geq 0
}
for all $x,y\in I$ with $x<y$, in other words, $f$ is $\omega$-Jensen convex on $I$ if and only if it is 
$\omega$-convex over all $(n+1)$-element arithmetic sequences in $I$. In the standard setting 
$\omega=(p_0,\dots,p_{n-1})$, these notions were introduced by Gil\'anyi--P\'ales in \cite{GilPal08}.

The following result shows that cyclic $(t,\omega)$-convexity implies $(r,\omega)$-convexity for all
$r\in\Q_+^n$. In the particular case $n=2$, and $\omega=(p_0,p_1)$, our result reduces to that of by Kuhn
\cite{Kuh84} and Dar\'oczy--P\'ales \cite{DarPal87}. For the higher-order case $\omega=(p_0,\dots,p_{n-1})$,
the analogous statement was established by Gil\'anyi--P\'ales \cite{GilPal08}.

\Thm{5}{Let $t\in\R_+^n$ and $f:I\to\R$. If $f$ is cyclically $(t,\omega)$-convex then it is
$(r,\omega)$-convex for all $r\in\Q_+^n$.}

\begin{proof}Let $r=(r_1,\dots,r_n)\in\Q_+^n$ and let $x\in I$ and $h>0$ with
$x+(r_1+\dots+r_n)h\in I$. There exist positive integers $q_1,\dots,q_n,q\in\N$ such that
\Eq{*}{
   r_1=\frac{q_1}{q},\qquad \dots, \qquad r_n=\frac{q_n}{q}.
}
Let us consider the elements
\Eq{*}{
  x_{kn+j}:=x + \bigg( k+\frac{1}{T}\sum_{i=1}^{j}t_i\bigg)\frac{h}{q}
}
for $k=0,\dots,(q_1+\dots+q_n)$, $j=0,\dots,k$, where $T:=\sum_{i=1}^{j}t_i$ and
we use the convention $\sum_{i=1}^{0}t_i:=0$.
It follows easily from this construction that
\Eq{*}{
  x_{i+1}-x_i=t_{\pi(1)}\frac{h}{Tq}, \qquad x_{i+2}-x_i=(t_{\pi(1)}+t_{\pi(2)})\frac{h}{Tq},
  \qquad \dots,\qquad x_{i+n}-x_i=(t_{\pi(1)}+\dots +t_{\pi(n)} )\frac{h}{Tq},
}
where $\pi$ is a cyclic permutation of $\{1,\dots,n\}$ for an arbitrary
$i\in \{0,\dots,\big((q_1+\dots+q_n)-1\big)n\}$. Therefore, the cyclic $(t,\omega)$-convexity
of $f$ implies that
\Eq{*}{
    \Phi_{(\omega,f)}(x_i,\dots,x_{i+n})\geq0 \qquad(i\in \{0,\dots,\big((q_1+\dots+q_n)-1\big)n\}).
}
Applying \cor{2}, it follows that
\Eq{*}{
    \Phi_{(\omega,f)}\big(x_0,x_{q_1n},\dots,x_{(q_1+\dots+q_n)n}\big)\geq0,
}
i.e.,
\Eq{*}{
    \Phi_{(\omega,f)}\big(x_0,x_0+r_1h,\dots,x_0+(r_1+\dots+r_n)h\big)\geq0,
}
which proves the $(r,\omega)$-convexity of $f$.
\end{proof}

The following corollary is an immediate consequence of \thm{5} and in the standard setting it was
established by Gilányi and Páles in \cite{GilPal08}.

\Cor{5}{If a function $f:I\to\R$ is $\omega$-Jensen convex, then it is $(r,\omega)$-convex for all
$r\in\Q_+^n$.}

Motivated by the above result, for a function $f:I\to\R$, we introduce the following notation
\Eq{*}{
  C_{(\omega,f)}:=\{t\in\R^n_+\mid f\mbox{ is cyclically $(t,\omega)$-convex}\}.
}
Obviously, $C_{(\omega,f)}$ is a cone, i.e., for $\lambda>0$, we have $\lambda C_{(\omega,f)}\subseteq
C_{(\omega,f)}$. The statement of \thm{5} says that the nonemptyness of $C_{(\omega,f)}$ implies
$\Q^n_+\subseteq C_{(\omega,f)}$. In the standard setting $n=2$, $\omega=(p_0,p_1)$, Kuhn's proved
(\cite{Kuh84}) that provided that $C_{(\omega,f)}\neq\emptyset$, there exists a subfield $\F$ of $\R$ such
that $C_{(\omega,f)}=\R_+\cdot\F_+^n$. It is an open problem if this statement remains valid in the general
case.

The next result shows that the $(t,\omega)$-convexity property is always the consequence of those
$(s,\omega)$-convexity properties where the coordinates of $s\in\R^n_+$ are equal to two subsequent
coordinates of $t$. To shorten the notation, we denote
\Eq{*}{
   (u\mathbf{1}_k,v\mathbf{1}_{m}):=(\underbrace{u,\dots,u}_{k\mbox{\rm\footnotesize-times}},
   \underbrace{v,\dots,v}_{m\mbox{\rm\footnotesize-times}}) \qquad(u,v\in\R,\,k,m\in\N).
}

\Thm{5+}{Assume that $n\geq2$ and let $(t_1,\dots,t_n)\in\R_+^n$. If, for all $(k,i)\in\{1,\dots,n-1\}^2$, 
the function $f:I\to\R$ is $\big((t_i\mathbf{1}_k,t_{i+1}\mathbf{1}_{n-k}),\omega\big)\mbox{-convex}$,
then $f$ is $\big((t_1,\dots,t_n),\omega\big)$-convex on $I$.}

\begin{proof}
We can assume that $n\geq3$, because if $n=2$, then the statement of the theorem is trivial. Let
$t_1,\dots,t_n>0$, $x\in I$ and $h>0$ with $x+(t_1+\cdots+t_n)h \in I$. Let us consider the elements
$x_0<x_1<\dots<x_{n(n-1)}$ defined by
\Eq{*}{
    x_{(n-1)(k-1)+j}:=x+\bigg((t_0+\dots+t_{k-1}) + \frac{j t_{k}}{n-1} \bigg)h
}
for all $(k,j)\in\big(\{1,\dots,n\}\times\{0,\dots,n-2\}\big)\cup\{(n,n-1)\}$, where $t_0:=0$. Notice that
for all elements $\ell\in\{0,\dots,n(n-1)-1\}$, there exists a unique pair
$(k_\ell,j_\ell)\in\{1,\dots,n\}\times\{0,\dots,n-2\}$ such that $\ell=(n-1)(k_\ell-1)+j_\ell$. On the other
hand, for $\ell\in\{0,\dots,n(n-2)\}$, we have
\Eq{*}{
    x_{\ell+1}-x_\ell&=\cdots=x_{(n-1)k_\ell}-x_{(n-1)k_\ell-1}&=\frac{t_{k_\ell}}{n-1},\\
    x_{(n-1)k_\ell+1}-x_{(n-1)k_\ell}&=\cdots=x_{\ell+n}-x_{\ell+n-1}&=\frac{t_{k_\ell+1}}{n-1}.
}
Therefore, using the 
$\big((t_{k_\ell}\mathbf{1}_{(n-1)k_\ell-\ell},t_{k_\ell+1}\mathbf{1}_{\ell+n-(n-1)k_\ell}), 
\omega\big)\big)$-convexity of $f$, we obtain that
\Eq{*}{
    \Phi_{(\omega,f)}(x_\ell,\dots,x_{\ell+n})\geq 0 \qquad (\ell\in\{0,\dots,n(n-2)\}).
}
Applying the \cor{2}, it follows that
\Eq{*}{
    \Phi_{(\omega,f)}(x_{0(n-1)},x_{1(n-1)},\dots,x_{n(n-1)}) \geq 0,
}
i.e.,
\Eq{*}{
    \Phi_{(\omega,f)}(x,x+t_1h,\dots,x+(t_1+\cdots+t_n)h)\geq 0,
}
which implies that the function $f$ is $\big((t_1,\dots,t_n),\omega\big)$-convex on $I$.
\end{proof}

The following two results are immediate consequences of \thm{5}.

\Cor{5+}{Assume that $n\geq2$ and let $(t_1,\dots,t_n)\in\R_+^n$. Denote $t_{n+1}:=t_1$. If, for all 
$(k,i)\in\{1,\dots,n-1\}\times\{1,\dots,n\}$, the function $f:I\to\R$ is 
$\big((t_i\mathbf{1}_k,t_{i+1}\mathbf{1}_{n-k}),\omega\big)\mbox{-convex}$,
then $f$ is cyclically $\big((t_1,\dots,t_n),\omega\big)$-convex on $I$.}

\Cor{5++}{Assume that $n\geq2$, and let $T\subseteq\R_+$ be a nonempty set. Then, for all
$t_1,\dots,t_n \in T$ the function $f:I\to \R$ is $\big((t_1,\dots,t_n),\omega\big)$-convex if and only if,
for all $t,s \in T$ and $k\in\{1,\dots,n\}$, the function $f$ is
$\big((t\mathbf{1}_k,s\mathbf{1}_{n-k}),\omega\big)\mbox{-convex}$.}

\section{Characterization of the convexity properties by Dinghas type derivatives}

In the rest of this paper, let $I\subset\R$ be a nondegenerate interval and let
$\omega=(\omega_1,\dots,\omega_n):I\to\R^n$ be an $n$-dimensional positive Chebyshev system over $I$ and 
choose $\omega_{n+1}:I\to\R$ such that $\omega^*=(\omega_1,\dots,\omega_n,\omega_{n+1})$ be an 
$(n+1)$-dimensional positive Chebyshev system over $I$.

For a point $p\in I$, define the following two lower Dinghas type generalized derivatives:
\Eq{*}{
 \underline{D}_{\,\omega^*}f(p)
     &:= \liminf_{\begin{smallmatrix}
          x_n-x_0\to 0 \\ (x_0,\dots,x_n)\in\sigma_{n+1}(I) \\ x_0\leq p\leq x_n
         \end{smallmatrix}
        }\big[x_0,\dots,x_n;f\big]_{\omega^*} \\
     &:= \lim_{\delta\to0^+}\inf\big\{\big[x_0,\dots,x_n;f\big]_{\omega^*}
          \,:\, (x_0,\dots,x_n)\in\sigma_{n+1}(I),\, x_0\leq p\leq x_n,\, x_n-x_0<\delta\big\}
}
and
\Eq{*}{
 \underline{D}_{\,(t,\omega^*)}f(p)
     &:= \liminf_{\begin{smallmatrix}
          y-x\to 0 \\ (x,y)\in\sigma_2(I)\\ x\leq p\leq y
         \end{smallmatrix}
        }\Big[x,\tfrac{t_2+\dots+t_n}{t_1+\dots+t_n}x+\tfrac{t_1}{t_1+\dots+t_n}y,
       \dots,\tfrac{t_n}{t_1+\dots+t_n}x+\tfrac{t_1+\dots+t_{n-1}}{t_1+\dots+t_n}y,y;f\Big]_{\omega^*} \\
     &:=
\lim_{\delta\to0^+}\inf\Big\{\Big[x,\tfrac{t_2+\dots+t_n}{t_1+\dots+t_n}x+\tfrac{t_1}{t_1+\dots+t_n}y,
       \dots,\tfrac{t_n}{t_1+\dots+t_n}x+\tfrac{t_1+\dots+t_{n-1}}{t_1+\dots+t_n}y,y;f\Big]_{\omega^*} \\
          &\hspace{6cm}\,:\, (x,y)\in\sigma_{2}(I),\, x\leq p\leq y,\, y-x<\delta\Big\},
}
where $(t_1,\dots,t_n)\in\R_+^n$. (For the original definition of A.\ Dinghas, we refer to the paper
\cite{Din66}.)

In the subsequent subsections we will obtain mean value inequalities and characterization theorems for
$\omega$-convexity, $\omega$-Jensen convexity, and for $((t_1,t_2),\omega)$-convexity, respectively, and we
will also establish that these convexity properties are localizable.

\subsection{Mean value inequality and characterization for $\omega$-convexity}

\Thm{6}{Let $f:I\to\R$. Then, for every $(x_0,\dots,x_n)\in\sigma_{n+1}(I)$, there exists $p\in[x_0,x_n]$
such that
\Eq{T6}{
    \big[x_0,\dots,x_n;f\big]_{\omega^*} \geq \underline{D}_{\,\omega^*}f(p).
}}

\begin{proof}
Let $(x_0,\dots,x_n)\in\sigma_{n+1}(I)$ be arbitrary. Denote $\max_{j\in\{0,\dots,n-1\}}(x_{j+1}-x_{j})$ by
$d$. We show the existence of a sequence $(x_0^{(k)},\dots,x_n^{(k)})\in\sigma_{n+1}(I)$ such that
$(x_0^{(1)},\dots,x_n^{(1)})=(x_0,\dots,x_n)$ and, for $k\in\N$,
\begin{eqnarray}
  &(x_0^{(k+1)},\dots,x_n^{(k+1)})\in
    \sigma_{n+1}\bigg(\Big\{x_0^{(k)},\tfrac{x_0^{(k)}+x_1^{(k)}}{2},x_1^{(k)},
     \tfrac{x_1^{(k)}+x_2^{(k)}}{2},\dots,\tfrac{x_{n-1}^{(k)}+x_n^{(k)}}{2},x_n^{(k)}\Big\}\bigg),
   \label{Ei1}\\[2mm]
  &\big[x_0^{(k+1)},\dots,x_n^{(k+1)};f\big]_{\omega^*}
     \leq \big[x_0^{(k)},\dots,x_n^{(k)};f\big]_{\omega^*}, \label{Ei2}\\[2mm]
  &\displaystyle\max_{i\in\{0,\dots,n-1\}} \Big(x_{i+1}^{(k+1)}-x_{i}^{(k+1)}\Big)\leq \frac{d}{2^k}.
   \label{Ei3}
\end{eqnarray}
Assume that we have already constructed $(x_0^{(k)},\dots,x_n^{(k)})\in\sigma_{n+1}(I)$. Define the points
$y_0^{(k)},\dots,y_{2n}^{(k)}$ by
\Eq{t1}{
  y_{2j}^{(k)}:=x_j^{(k)} \quad (j\in\{0,\dots,n\}) \qquad\mbox{and}\qquad
  y_{2j-1}^{(k)}:=\tfrac{x_{j-1}^{(k)}+x_j^{(k)}}{2} \quad (j\in\{1,\dots,n\}).
}
Then, obviously, $y_{0}^{(k)}<y_{1}^{(k)}<y_{2}^{(k)}<\dots<y_{2n}^{(k)}$ and
\Eq{t2}{
  \max_{j\in\{0,\dots,2n-1\}} \Big(y_{j+1}^{(k)}-y_{j}^{(k)}\Big)
  \leq \frac12\cdot \max_{i\in\{0,\dots,n-1\}} \Big(x_{i+1}^{(k)}-x_{i}^{(k)}\Big)\leq \frac{d}{2^k}.
}
By the Generalized Chain Inequality, i.e., by \cor{1},
\Eq{t3}{
  \big[x_0^{(k)},\dots,x_n^{(k)};f\big]_{\omega^*} \geq
  \min_{i\in\{0,\dots,n\}}\big[y_{i}^{(k)},\dots,y_{i+n}^{(k)};f\big]_{\omega^*}.
}
Now, choose $j\in\{0,\dots,n\}$ such that
$\big[y_{j}^{(k)},\dots,y_{j+n}^{(k)};f\big]_{\omega^*}=\min_{i\in\{0,\dots,n\}}\big[y_{i}^{(k)},\dots,y_{i+n}
^{(k)};f\big]_{\omega^*}$ and define $(x_0^{(k+1)},\dots,x_n^{(k+1)}):=(y_{j}^{(k)},\dots,y_{j+n}^{(k)})$.
Then, in view of properties \eq{t1}, \eq{t3}, and \eq{t2}, we obtain that \eq{i1}, \eq{i2}, and \eq{i3} hold,
respectively.

Observe that $\big(x_0^{(k)}\big)$ is a nondecreasing and $\big(x_n^{(k)}\big)$ is a nonincreasing sequence,
furthermore,
\Eq{*}{
  0<x_n^{(k)}-x_0^{(k)}
  <n \cdot \max_{i\in\{0,\dots,n-1\}} \Big(x_{i+1}^{(k)}-x_{i}^{(k)}\Big)\leq\frac{nd}{2^{k-1}}.
}
Hence the sequences $\big(x_0^{(k)}\big)$ and $\big(x_n^{(k)}\big)$ have a common limit point which we denote
by $p$. Using \eq{i2},
\Eq{*}{
  \big[x_0,\dots,x_n;f\big]_{\omega^*}
  =\big[x_0^{(1)},\dots,x_n^{(1)};f\big]_{\omega^*}
  \geq \liminf_{k\to\infty}\big[x_0^{(k)},\dots,x_n^{(k)};f\big]_{\omega^*}
  \geq \underline{D}_{\,\omega^*}f(p),
}
which results \eq{T6}.
\end{proof}

The next result characterizes the $\omega$-convexity property in terms of the nonnegativity of the
generalized lower Dinghas type derivative $\underline{D}_{\,\omega^*}$. 


\Cor{6}{A function $f:I\to\R$ is $\omega$-convex on $I$ if and only if, for all $p\in I$,
\Eq{Df}{
  \underline{D}_{\,\omega^*}f(p)\geq0.
}}

\begin{proof} Assume that $f$ is $\omega$-convex. Then, for all $(x_0,\dots,x_n)\in\sigma_{n+1}(I)$,
we have that $\big[x_0,\dots,x_n;f\big]_{\omega^*}\geq 0$, which, by the definition of the derivative
$\underline{D}_{\,\omega^*}f(p)$, implies that \eq{Df} holds for all $p\in I$.

To prove the reversed statement, assume that \eq{Df} holds. Choose $(x_0,\dots,x_n)\in\sigma_{n+1}(I)$
arbitrarily. Then, by \thm{6}, there exists $p\in I$ such that \eq{T6} holds which, together with \eq{Df},
implies $\big[x_0,\dots,x_n;f\big]_{\omega^*}\geq 0$. This is equivalent to the inequality
$\Phi_{(\omega,f)}(x_0,\dots,x_n)\geq 0$. This proves the $\omega$-convexity of $f$.
\end{proof}

The next immediate consequence shows that $\omega$-convexity is a localizable property.

\Cor{7}{A function $f:I\to\R$ is $\omega$-convex on $I$ if and only if, for all $p\in I$,
there exists a neighborhood $U$ of $p$ such that $f$ is $\omega$-convex on $U$.}

\subsection{Mean value inequality and characterization for $\omega$-Jensen convexity}

The following result is analogous to \thm{6} and establishes a mean value inequality for $\omega$-Jensen type
divided differences.

\Thm{7}{Let $f:I\to\R$. Then, for every $x,y\in I$ with $x<y$, there exists $p\in[x,y]$ such that
\Eq{T7}{
\big[x,\tfrac{n-1}{n}x+\tfrac{1}{n}y,\dots,\tfrac{1}{n}x+\tfrac{n-1}{n}y,y;f\big]_{\omega^*} \geq
\underline{D}_{\,(\mathbf{1}_n,\omega^*)}f(p).
}}

\begin{proof}
Let $x,y\in I$ be arbitrary with $x<y$.
We prove by induction that there exists a sequence $\big(x^{(k)}\big)$ in $I$ such that $x^{(1)}=x$ and, for $k\in\N$,
\Eq{xk}{
  x^{(k+1)}&\in
    \Big\{x^{(k)},x^{(k)}+\tfrac{1}{n2^k}(y-x),\dots,x^{(k)}+\tfrac{n}{n2^k}(y-x)\Big\},\\[2mm]
  \big[x^{(k+1)},x^{(k+1)}+\tfrac{1}{n2^k}(y-x),&\dots,x^{(k+1)}+\tfrac{n}{n2^k}(y-x);f\big]_{\omega^*}\\
     &\leq
   \big[x^{(k)},x^{(k)}+\tfrac{1}{n2^{k-1}}(y-x),\dots,x^{(k)}+\tfrac{n}{n2^{k-1}}(y-x);f\big]_{\omega^*},
}
Assume that we have already constructed $x^{(k)}\in I$. Define the points $y_0^{(k)},\dots,y_{2n}^{(k)}$ by
\Eq{tt1}{
  y_{j}^{(k)}:=x^{(k)}+\tfrac{j}{n2^k}(y-x) \qquad (j\in\{0,\dots,2n\}).
}
Then, obviously, $y_{0}^{(k)}<y_{1}^{(k)}<\dots<y_{2n}^{(k)}$ and
\Eq{*}{
  \big(x^{(k)},x^{(k)}+\tfrac{1}{n2^{k-1}}(y-x),\dots,x^{(k)}+\tfrac{n}{n2^{k-1}}(y-x)\big)
  \in\sigma_{n+1}\big(\big\{y_{0}^{(k)},y_{1}^{(k)},\dots,y_{2n}^{(k)}\big\}\big).
}
Therefore, by applying the Generalized Chain Inequality, it follows that
\Eq{tt3}{
  \big[x^{(k)},x^{(k)}+\tfrac{1}{n2^{k-1}}(y-x),\dots,x^{(k)}+\tfrac{n}{n2^{k-1}}(y-x);f\big]_{\omega^*}
  \geq \min_{i\in\{0,\dots,n\}}\big[y_{i}^{(k)},\dots,y_{i+n}^{(k)};f\big]_{\omega^*}.
}
Now, choose $j\in\{0,\dots,n\}$ such that
$\big[y_{j}^{(k)},\dots,y_{j+n}^{(k)};f\big]_{\omega^*}=\min_{i\in\{0,\dots,n\}}\big[y_{i}^{(k)},\dots,y_{i+n}
^{(k)};f\big]_{\omega^*}$ and define $x^{(k+1)}:=y_{j}^{(k)}$.
Then, the properties \eq{tt1} and \eq{tt3} imply that the two conditions in \eq{xk} hold.

Observe that $\big(x^{(k)}\big)$ is a nondecreasing and
$\big(y^{(k)}\big):=\big(x^{(k)}+\tfrac{n}{n2^{k-1}}(y-x)\big)$ is a nonincreasing sequence, furthermore,
\Eq{*}{
  0<y^{(k)}-x^{(k)}=\tfrac{y-x}{2^{k-1}}.
}
Thus these sequences have a common limit point which we denote by $p$. Using the second inequality in
\eq{xk},
\Eq{*}{
  &\big[x,\tfrac{n-1}{n}x+\tfrac{1}{n}y,\dots,\tfrac{1}{n}x+\tfrac{n-1}{n}y,y;f\big]_{\omega^*}
   =\big[x^{(1)},x^{(1)}+\tfrac{1}{n}(y-x),\dots,x^{(1)}+\tfrac{n}{n}(y-x);f\big]_{\omega^*} \\
  &\geq\liminf_{k\to\infty}\big[x^{(k)},x^{(k)}+\tfrac{1}{n2^{k-1}}(y-x),\dots,
                           x^{(k)}+\tfrac{n}{n2^{k-1}}(y-x);f\big]_{\omega^*}  \\
  &= \liminf_{k\to\infty}\big[x^{(k)},\tfrac{n-1}{n}x^{(k)}+\tfrac{1}{n}y^{(k)},\dots,
	  \tfrac{1}{n}x^{(k)}+\tfrac{n-1}{n}y^{(k)},y^{(k)};f\big]_{\omega^*}
  \geq \underline{D}_{\,(\mathbf{1},\omega^*)}f(p),
}
which results \eq{T7}.
\end{proof}

The following two corollaries directly follow from \thm{7} as \cor{6} and \cor{7} from \thm{6}.

\Cor{8}{A function $f:I\to\R$ is $\omega$-Jensen convex on $I$ if and only if, for all $p\in I$,
\Eq{*}{
  \underline{D}_{\,(\mathbf{1}_n,\omega^*)}f(p)\geq0.
}}

\Cor{9}{A function $f:I\to\R$ is $\omega$-Jensen convex on $I$ if and only if, for all $p\in I$, there
exists a neighborhood $U$ of $p$ such that $f$ is $\omega$-Jensen convex on $U$.}

\subsection{Mean value inequality and characterization for $((t_1,t_2),\omega)$-convexity}

\Thm{9}{Assume that $n=2$ and let $(t_1,t_2)\in\R_+^2$ and $f:I\to\R$. Then, for every $(x,y)\in \sigma_2(I)$,
there exists $p\in[x,y]$ such that
\Eq{T9}{
\big[x,\tfrac{t_2}{t_1+t_2}x+\tfrac{t_1}{t_1+t_2}y,y;f\big]_{\omega^*} \geq
\min \big\{\underline{D}_{\,((t_1,t_2),\omega^*)}f(p),\underline{D}_{\,((t_2,t_1),\omega^*)}f(p)\big\}.
}}

\begin{proof}
Let $t_1,t_2\in \R_+$ and $(x,y)\in \sigma_2(I)$ be arbitrary. We may assume that $t_1<t_2$. We prove by
induction that there exists a sequence $ \big(x_0^{(k)},x_1^{(k)},x_2^{(k)}\big)\in\sigma_3(I) $ such that
\Eq{*}{
\big(x_0^{(1)},x_1^{(1)},x_2^{(1)}\big)=\big(x,\tfrac{t_2}{t_1+t_2}x+\tfrac{t_1}{t_1+t_2}y,y\big)
}
and, for all $k\in \N$,
\begin{eqnarray}
    &\big(x_0^{(k+1)},x_1^{(k+1)},x_2^{(k+1)}\big) \in
	\sigma_3\Big(\Big\{ x_0^{(k)},\tfrac{t_2}{t_1+t_2}x_0^{(k)}+\tfrac{t_1}{t_1+t_2}x_1^{(k)},
	\tfrac{t_1}{t_1+t_2}x_0^{(k)}+\tfrac{t_2}{t_1+t_2}x_1^{(k)},x_1^{(k)}, \nonumber\\
	&  \hspace{5cm}\tfrac{t_2}{t_1+t_2}x_1^{(k)}+\tfrac{t_1}{t_1+t_2}x_2^{(k)},
	  \tfrac{t_1}{t_1+t_2}x_1^{(k)}+\tfrac{t_2}{t_1+t_2}x_2^{(k)},x_2^{(k)}\Big\}\Big),  \label{EE9i1}\\
    &\frac{x_1^{(k+1)}-x_0^{(k+1)}}{x_2^{(k+1)}-x_1^{(k+1)}}\in\big\{\frac{t_1}{t_2},\frac{t_2}{t_1}\big\},
    \label{EE9i1+}\\[2mm]
    &\big[x_0^{(k+1)},x_1^{(k+1)},x_2^{(k+1)};f\big]_{\omega^*}\leq \big[x_0^{(k)},x_1^{(k)},x_2^{(k)};f\big]_{\omega^*},   \label{EE9i2}\\[2mm]
    &x_2^{(k+1)}-x_0^{(k+1)} \leq
        \big(\frac{t_2}{t_1+t_2}\big)^k \cdot (y-x). \label{EE9i3}
\end{eqnarray}
Observe that \eq{E9i1+} and \eq{E9i3} hold for $k=0$.
Assume that we have already constructed $ \big(x_0^{(k)},x_1^{(k)},x_2^{(k)}\big)\in\sigma_3(I)$.
By \eq{E9i1+}, we have that $\frac{x_1^{(k)}-x_0^{(k)}}{x_2^{(k)}-x_1^{(k)}}=\frac{t_i}{t_j}
$ for some $(i,j)\in\{(1,2),(2,1)\}$. This means that
\Eq{ij}{
   x_1^{(k)}=\tfrac{t_j}{t_1+t_2}x_0^{(k)}+\tfrac{t_i}{t_1+t_2}x_2^{(k)}.
}
Now define the points $y_0^{(k)},\dots,y_4^{(k)} $ by
\Eq{t91}{
  &y_{2\ell}^{(k)}:=x_\ell^{(k)} \qquad (\ell\in\{0,1,2\}),\\
  &y_{1}^{(k)}:=\tfrac{t_j}{t_1+t_2}x_0^{(k)}+\tfrac{t_i}{t_1+t_2}x_1^{(k)},\\
  &y_{3}^{(k)}:=\tfrac{t_i}{t_1+t_2}x_1^{(k)}+\tfrac{t_j}{t_1+t_2}x_2^{(k)}.
}
Then, obviously $y_0^{(k)}<\dots<y_4^{(k)}$. Using \eq{ij} and then \eq{E9i3}, for
$\ell\in\{0,1,2\}$, we get
\Eq{t92}{
  y_{\ell+2}^{(k)}-y_\ell^{(k)}
  &\leq\max\big\{y_2^{(k)}-y_0^{(k)},y_3^{(k)}-y_1^{(k)},y_4^{(k)}-y_2^{(k)}\big\}\\
  &=\max\big\{x_1^{(k)}-x_0^{(k)},\tfrac{t_i}{t_1+t_2}x_1^{(k)}+\tfrac{t_j}{t_1+t_2}x_2^{(k)}
     -\tfrac{t_j}{t_1+t_2}x_0^{(k)}-\tfrac{t_i}{t_1+t_2}x_1^{(k)},x_2^{(k)}-x_1^{(k)}\big\}\\
  &=\max\big\{\tfrac{t_i}{t_1+t_2},\tfrac{t_j}{t_1+t_2},\tfrac{t_j}{t_1+t_2}\big\}
   \big(x_2^{(k)}-x_0^{(k)}\big)\\
  &=\tfrac{t_2}{t_1+t_2}\big(x_2^{(k)}-x_0^{(k)}\big)
  \leq \big(\tfrac{t_2}{t_1+t_2}\big)^k(y-x).
}
Furthermore, using \eq{t91} and \eq{ij}, we also have that
\Eq{t94}{
  \frac{y_{1}^{(k)}-y_{0}^{(k)}}{y_{2}^{(k)}-y_{1}^{(k)}}=\frac{t_i}{t_j}, \qquad
  \frac{y_{2}^{(k)}-y_{1}^{(k)}}{y_{3}^{(k)}-y_{2}^{(k)}}=\frac{t_i}{t_j}, \qquad
  \frac{y_{3}^{(k)}-y_{2}^{(k)}}{y_{4}^{(k)}-y_{3}^{(k)}}=\frac{t_j}{t_i}.
}
By the Generalized Chain Inequality, i.e., by \cor{1},
\Eq{t93}{
  \big[x_0^{(k)},x_1^{(k)},x_2^{(k)};f\big]_{\omega^*} \geq
  \min_{\ell\in\{0,1,2\}}\big[y_{\ell}^{(k)},y_{\ell+1}^{(k)},y_{\ell+2}^{(k)};f\big]_{\omega^*}.
}
Choose $m\in\{0,1,2\}$ such that $\big[y_{m}^{(k)},y_{m+1}^{(k)},y_{m+2}^{(k)};f\big]_{\omega^*}
=\min_{\ell\in\{0,1,2\}}\big[y_{\ell}^{(k)},y_{\ell+1}^{(k)},y_{\ell+2}^{(k)};f\big]_{\omega^*}$ and define
$(x_0^{(k+1)},x_1^{(k+1)},x_2^{(k+1)}):=(y_{m}^{(k)},y_{m+1}^{(k)},y_{m+2}^{(k)})$. Then, in view of
properties \eq{t91}, \eq{t94}, \eq{t93}, and \eq{t92}, we obtain that \eq{E9i1}, \eq{E9i1+}, \eq{E9i2}, and
\eq{E9i3} hold,
respectively.

Observe that $\big(x_0^{(k)}\big)$ is a nondecreasing and $\big(x_2^{(k)}\big)$ is a nonincreasing sequence,
furthermore,
\Eq{*}{
  0<x_2^{(k)}-x_0^{(k)}
  \leq \big(\tfrac{t_2}{t_1+t_2}\big)^{k-1}\cdot (y-x).
}
Hence the sequences $\big(x_0^{(k)}\big)$ and $\big(x_2^{(k)}\big)$ have a common limit point which we denote
by $p$. By property \eq{E9i1+}, one of the sets
\Eq{*}{
  K_1:=\bigg\{k\in\N:\frac{x_1^{(k)}-x_0^{(k)}}{x_2^{(k)}-x_1^{(k)}}=\frac{t_1}{t_2}\bigg\},\qquad
  K_2:=\bigg\{k\in\N:\frac{x_1^{(k)}-x_0^{(k)}}{x_2^{(k)}-x_1^{(k)}}=\frac{t_2}{t_1}\bigg\}
}
is infinite. Assume that $K_1$ is infinite let $k_1<k_2<\cdots<k_n<\cdots$
be chosen such that $K_1=\{k_1,k_2,\dots,k_n,\dots\}$. Using \eq{E9i2},
\Eq{*}{
  \big[x,\tfrac{t_2}{t_1+t_2}x+\tfrac{t_1}{t_1+t_2}y,y;f\big]_{\omega^*}
  &=\big[x_0^{(1)},x_1^{(1)},x_2^{(1)};f\big]_{\omega^*}
  \geq \liminf_{n\to\infty}\big[x_0^{(k_n)},x_1^{(k_n)},x_2^{(k_n)};f\big]_{\omega^*}\\
  &= \liminf_{n\to\infty}\big[x_0^{(k_n)},
    \tfrac{t_2}{t_1+t_2}x_0^{(k_n)}+\tfrac{t_1}{t_1+t_2}x_2^{(k_n)},x_2^{(k_n)};f\big]_{\omega^*}
  \geq \underline{D}_{\,((t_1,t_2),\omega^*)}f(p),
}
If $K_2$ is infinite, then in a similar manner, we get that
\Eq{*}{
  \big[x,\tfrac{t_2}{t_1+t_2}x+\tfrac{t_1}{t_1+t_2}y,y;f\big]_{\omega^*}
  \geq \underline{D}_{\,((t_2,t_1),\omega^*)}f(p),
}
which results \eq{T9}.
\end{proof}

The following two corollaries directly follow from \thm{9} as \cor{6} and \cor{7} from \thm{6}.
The first corollary characterizes symmetric $((t_1,t_2),\omega)$-convexity in terms of the nonnegativity of
two generalized derivatives. The second one shows that symmetric $((t_1,t_2),\omega)$-convexity is a
localizable property.

\Cor{10}{Assume that $n=2$ and let $(t_1,t_2)\in\R_+^2$. Then a function $f:I\to\R$ is
symmetrically $((t_1,t_2),\omega)$-convex on $I$ if and only if, for all $p\in I$,
\Eq{*}{
  \underline{D}_{\,((t_1,t_2),\omega^*)}f(p)\geq0 \qquad\mbox{and}\qquad
  \underline{D}_{\,((t_2,t_1),\omega^*)}f(p)\geq0.
}}

\Cor{11}{Assume that $n=2$ and let $(t_1,t_2)\in\R_+^2$. Then a function $f:I\to\R$ is symmetrically
$((t_1,t_2),\omega) $-convex on $I$ if and only if, for all $p\in I$, there exists a neighborhood $U$
of $p$ such that $f$ is symmetrically $((t_1,t_2),\omega) $-convex on $U$.}


\providecommand{\bysame}{\leavevmode\hbox to3em{\hrulefill}\thinspace}
\providecommand{\MR}{\relax\ifhmode\unskip\space\fi MR }
\providecommand{\MRhref}[2]{%
  \href{http://www.ams.org/mathscinet-getitem?mr=#1}{#2}
}
\providecommand{\href}[2]{#2}

\end{document}